\documentclass[11pt]{article}
\usepackage{geometry}
\usepackage{amsmath,amsfonts,amsthm,amssymb}
\title{Scattering and modified scattering for abstract wave equations
with time-dependent dissipation}
\author{Jens Wirth\thanks{Institute of Applied Analysis, 
TU Bergakademie Freiberg, 09596 Freiberg, Germany}\ \thanks{\tt 
wirth@math.tu-freiberg.de}}
\date{}
\newtheorem{thm}{Theorem}[section]
\newtheorem{lem}[thm]{Lemma}
\newtheorem{cor}[thm]{Corollary}
\theoremstyle{remark}
\newtheorem{rem}{Remark} 
\newtheorem{expl}{Example} 
\DeclareMathOperator{\Ker}{Ker}
\DeclareMathOperator{\cl}{cl}
\def\Dom{\mathcal D}
\def\d{\mathrm d}
\def\R{\mathbb R}
\numberwithin{equation}{section}
\begin{document}
\maketitle

\begin{abstract}
  We consider the initial-value problem of abstract wave equations with
  weak dissipation. 

  We show that under conditions on the dissipation coefficient and its
  derivative the solutions to the abstract dissipative equation are closely
  related to solutions of the free problem multiplied by a decay function. 
  This paper gives the counterpart 
  to a recent paper of T.Yamazaki [Adv. Differential Equ., 11(4):419--456, 2006],
  where effective dissipation terms and the relation to the corresponding 
  abstract parabolic problem are considered. 

\end{abstract}

\section{Introduction }
\label{sec1}
Let $H$ be a separable Hilbert space with norm $\|\cdot\|$. Let further $A$ be a 
closed and self-adjoint non-negative operator with domain $\Dom(A)\subseteq H$. Note,
that the domain of $A$ becomes itself a Hilbert space if we endow it with the
graph norm $\|u\|_A^2=\|u\|^2+\|Au\|^2\sim\|\langle A\rangle u\|^2$. To simplify notation we 
set for the following $\Lambda=\sqrt{A}$.

We consider for a positive valued $C^1$-function $b=b(t)$ on $[0,\infty)$ the
abstract dissipative wave equation
\begin{equation}\label{eq:CP}
  \begin{cases}
    u''+2b(t) u'+ Au=0\\u(0)=u_1\in\Dom(\Lambda),\; u'(0)=u_2\in H
  \end{cases}
\end{equation}
and compare its solutions to the corresponding free problem
\begin{equation}\label{eq:CPfree}
  \begin{cases}
    v''+ Av=0\\v(0)=v_1\in\Dom(\Lambda),\; v'(0)=v_2\in H
  \end{cases}
\end{equation}
in the abstract energy space associated to $\Lambda$ under decay assumptions on 
the coefficient function $b$. 

The first result is as follows. If we assume that $b\in L^1[0,\infty)$, then the 
solutions of \eqref{eq:CP} are asymptotically free in the energy norm. This 
means, if we set $E$ to be the closure of $(\Dom(\Lambda)\ominus{\Ker\Lambda})\times H$ with respect to the 
norm
\begin{equation}
  \|(u_1,u_2)\|_E = \|(\Lambda u_1,u_2)\|_{H\times H},
\end{equation}
the following theorem is valid:
\begin{thm}\label{thm1}
  Assume $b\in L^1[0,\infty)$. Then there exists an invertible operator 
  $W_+\in\mathcal L(E)$ of the abstract energy space, such that for
  $(u_1,u_2)\in E$ and $(v_1,v_2)= W_+(u_1,u_2)$ the corresponding solutions to 
  \eqref{eq:CP} and \eqref{eq:CPfree} satisfy
  \begin{equation}\label{eq:1.4}
    \|(u,u')-(v,v')\|_E\leq C\|(u_1,u_2)\|_E \int_t^\infty b(\tau)\d\tau\to 0
  \end{equation}
  as $t\to\infty$. 
\end{thm}
\begin{rem}
  Note that the solutions to the free problem have the conservation of energy property, i.e. it holds
  for all $t\in\R$ that $\|(v,v')\|_E=\|(v_1,v_2)\|_E$. Thus the statement of Theorem~\ref{thm1} implies the non-decay to zero
  of the energy of solutions to the perturbed problem \eqref{eq:CP} together with the existence of 
  asymptotically equivalent solutions to both problems.
\end{rem}
\begin{rem}
  Note, that the operator $W_+$ in this theorem is the inverse of the M\o{}ller wave operator. In this point our notation differs from \cite{LP67}.  This inverse operator is more appropriate to describe asymptotic properties of solutions to problems with perturbed coefficients.
\end{rem}

The aim of this article is to generalise this result to the case of so-called
non-effective weak dissipation terms $b=b(t)$ with
\begin{description}
\item[(A1)] $b(t)\geq0$ for all $t\in[0,\infty)$,
\item[(A2)] $|b(t)|\leq C_1 \langle t\rangle^{-1}$ and $|b'(t)|\leq C_2 \langle t\rangle^{-2}$,
\end{description}
including non-integrable coefficients $b\not\in L^1[0,\infty)$. In this case
the abstract energy $\|(u,u')\|_E$ of the solution decays to zero. We will show that if $0$ is 
not an eigenvalue of $A$ then the solutions can be written asymptotically as product of a free 
solution and a time-dependent function describing the energy decay rate. To be more precise,
let us formulate a theorem (which is contained in Theorem~\ref{thm:3.1}).
\begin{thm}\label{thm2}
  Assume (A1), (A2) together with $\limsup_{t\to\infty}tb(t)<\frac12$ and $\Ker A=\{0\}$. Then to 
  $u_1\in\Dom(\Lambda)$ and $u_2\in H$ there exist corresponding data $(v_1,v_2)=W_+(u_1,u_2)\in E$
  such that the corresponding solutions to \eqref{eq:CP} and \eqref{eq:CPfree} satisfy
  \begin{equation}
    \|\lambda(t)(u,u')-(v,v')\|_E\to0
  \end{equation}
  as $t\to\infty$, where the auxiliary function $\lambda(t)$ is given by $\lambda(t)=\exp(\int_0^tb(\tau)\d\tau)$.
\end{thm}
\begin{rem}
  If $b\not\in L^1[0,\infty)$ the integral in equation~\eqref{eq:1.4} becomes infinite. In order
  to obtain a relation between the solutions $u$ and $v$ we have to pay. This can be seen from
  the differences of the statements. The operator $W_+$ maps $\Dom(\Lambda)\times H$ continuously to $E$, it is in general not
  bounded from $E$ to $E$. Furthermore, we cannot give a uniform rate of asymptotic equivalence.
  The auxiliary function $\lambda(t)$ describes the decay of energy. If $\Ker A=\{0\}$, we get the two-sided estimate
  \begin{equation}
    \|(u,u')\|_E\sim\frac1{\lambda(t)}.
  \end{equation}
  Since
  the convergence rate in Theorem~\ref{thm2} is not uniform in the data, the constants in this estimate
  depend in a (non-linear) way on the data $(u_1,u_2)$ and do not yield  a two-sided norm estimate
  for a corresponding solution operator.
\end{rem}
\begin{expl}
  One typical class of coefficient functions $b=b(t)$ for our approach is 
  \begin{equation}
    b(t)=\frac\mu{(1+t)\log(e+t)\cdots\log^{[n]}(e^{[n]}+t)},\qquad \mu\geq0,\quad n=0,1,\ldots
  \end{equation}
  with iterated logarithms $\log^{[1]}=\log$, $\log^{[n+1]}=\log\circ\log^{[n]}$ and
  corresponding iterated exponentials $e^{[0]}=1$, $e^{[n+1]}=e^{(e^{[n]})}$. In this case
  $\lambda(t)\sim\left(\log^{[n]}(e^{[n]}+t\right)^\mu$. The decay rate for the abstract energy
  becomes arbitrary small in the scale of iterated logarithms.\hfill$\square$
\end{expl}

Recently, T.Yamazaki, \cite{Yam06}, \cite{Yam06b}, considered the related problem in the case of
effective weak dissipation, i.e. if $tb(t)\to\infty$ as $t\to\infty$ together with $b(t)\to0$ and suitable conditions on the first derivative. 
In this case, a relation to the free problem \eqref{eq:CPfree} does not occur. Instead ,
there arises a close relation to the abstract parabolic equation $b(t)w'+Aw=0$. While the purpose of
\cite{Yam06}, \cite{Yam06b} was to demonstrate the close relation of abstract wave equations with effective
dissipation and the corresponding abstract parabolic problem, our purpose is to show what happens if the
influence of the dissipation becomes less strong and how the abstract parabolic type asymptotics changes
to an abstract wave type asymptotics. The used methods are based on \cite{CH03}.

This paper is organised as follows. First we will sketch one proof of
Theorem~\ref{thm1} and give some remarks on possible generalisations. Then
in Section~\ref{sec2} we will summarise the construction of representations of
solutions to \eqref{eq:CP} by means of a spectral resolution of the operator 
$A$ and discuss its consequences for the decay rates of the energy. In 
Section~\ref{sec3} we state the main result of this paper and in 
Section~\ref{sec4} we discuss some possible applications of the abstract 
results in concrete settings.

The main part of this paper follows the author's treatment from \cite{Wir06} combined with ideas from \cite{CH03}. In that paper precise $L^p$--$L^q$ decay estimates are proven for solutions to the Cauchy problem of a wave equation with non-effective time-dependent dissipation. The approach was based on a diagonalisation scheme applied to the full symbol of the operator. Here we formulate the Cauchy problem in an abstract setting and base our representation of solutions on a spectral theorem for the operator $A$.  Besides this, there are two essential differences to the treatment in \cite{Wir06}. On the one hand the derivation of dispersive estimates is based on stronger assumptions on the coefficient function (we need estimates for $n+2$ derivatives of the coefficient function for dispersive estimates) and requires more steps of diagonalisation, while in the Hilbert space setting estimates for $b$ and its first derivative are sufficient. Second difference to \cite{Wir06} is that we obtain the modified scattering result not only for $\limsup_{t\to\infty}tb(t)< 1/2$. We use a new idea to impose conditions on the data in order to obtain a similar result for $\liminf_{t\to\infty}tb(t)>1/2$.

\subsection{Outline of the proof of the classical scattering result}

Theorem~\ref{thm1} can be reduced to abstract results already known from the literature. Using a substitution 
of the time variable $t=t(\tau)=\int_0^\tau a(s)\d s$ such that $\log \lambda(t)= \int_0^t b(s)\d s=\log a(\tau)$ reduces \eqref{eq:CP} 
to the abstract problem $w''+a^2(\tau)Aw=0$. Because we have $a^2(\tau)\to\int_0^\infty b(t)\d t<\infty$
and $\int_0^\infty (a(\tau)-a(\infty))\d\tau=\int_0^\infty \tau(t) b(t)\d t<\infty$, the technique of A.Arosio, \cite{Aro84}, can be applied. Transformation
back yields Theorem~\ref{thm1}.

For completeness and because the idea of the proof is basic for our treatment in sections \ref{sec2} and \ref{sec3}
we give an independent proof of Theorem~\ref{thm1}.

{\em Proof of Theorem~\ref{thm1}.} We use the canonic identification $J:E\stackrel\simeq\longrightarrow (H\ominus\Ker\Lambda)\times H$ of  the energy space with $\cl R(\Lambda)\times H$ and write the Cauchy problems in system form.
This yields
\begin{equation}\label{eq:system}
  \frac{\d}{\d t}  \begin{pmatrix}\Lambda u\\u'\end{pmatrix}
  =  \begin{pmatrix}\Lambda u'\\u''\end{pmatrix}
  =  \begin{pmatrix}&\Lambda\\-\Lambda &-2b(t)\end{pmatrix} 
  \begin{pmatrix}\Lambda u\\u'\end{pmatrix}.
\end{equation}
We denote by $\mathcal E(t,s)$ and $\mathcal E_0(t-s)$ the corresponding 
semi-groups of operators,
\begin{align}
 \frac{\d}{\d t} \mathcal E_0(t) 
 &=  \begin{pmatrix}&\Lambda\\-\Lambda &\end{pmatrix}  \mathcal E_0(t),
 \qquad &\mathcal E_0(0)&=I\in\mathcal L(H\times H), \\
 \frac{\d}{\d t} \mathcal E(t,s) 
 &=  \begin{pmatrix}&\Lambda\\-\Lambda &-2b(t)\end{pmatrix}  \mathcal E(t,s), 
\qquad &\mathcal E(s, s)&=I\in\mathcal L(H\times H) . 
\end{align}
The free semi-group $\mathcal E_0(t)$ is unitary (because its generator is
skew); while under assumption (A1) the semi-group $\mathcal E(t,s)$ is
contractive.\footnote{For a treatment without (A1) we construct 
$\mathcal E(t,s)=\mathcal E_0(t-s)\mathcal Q(t,s)$ directly by the aid of
formula \eqref{eq:Qdef}} If we consider the operator
\begin{equation}
  \mathcal Q(t,s)=\mathcal E_0(s-t) \mathcal E(t,s)\in\mathcal L(H\times H),
\end{equation}
we obtain the abstract differential equation
\begin{align}
  \frac{\d}{\d t}\mathcal Q(t,s)
  &=-\mathcal E_0(s-t)\begin{pmatrix}&\Lambda\\-\Lambda &\end{pmatrix}
  \mathcal E(t,s)+\mathcal E_0(s-t)\begin{pmatrix}&\Lambda\\-\Lambda &-2b(t)\end{pmatrix}
  \mathcal E(t,s)\notag\\
  &=-2b(t)\mathcal E_0(s-t)\begin{pmatrix}0&0\\0&1\end{pmatrix}\mathcal E_0(t-s)
  \mathcal Q(t,s)=\mathcal R(t,s)\mathcal Q(t,s)
\end{align}
with initial data $\mathcal Q(s, s)=I\in\mathcal L(H\times H)$. Using that the 
operator norm satisfies $|||\mathcal R(t,s)|||\leq2b(t)\in L^1[0,\infty)$ we see that
$\mathcal Q(t,s)$ can be represented by the abstract Peano-Baker formula 
\begin{equation}
  \label{eq:Qdef}
  \mathcal Q(t,s)=I+\sum_{k=1}^\infty \int_s^t \mathcal R(t_1,s) \int_s^{t_1} \mathcal R(t_2,s)\ldots
  \int_s^{t_{k-1}} \mathcal R(t_k,s)\d t_k\ldots\d t_1
\end{equation}
in terms of Bochner integrals. The $L^1$-bound on the norm of $\mathcal R$ 
implies a uniform bound on $\mathcal Q(t,s)$ and that
\begin{equation}
  \lim_{t\to\infty}\mathcal Q(t,s)=\mathcal Q(\infty,s)
\end{equation}
exists in $\mathcal L(H\times H)$ and satisfies
\begin{align}
  |||\mathcal Q(\infty,s)-\mathcal Q(t,s)|||\leq 2\int_t^\infty b(\tau)\d\tau\;\exp\left(2\int_0^\infty b(\theta)\d\theta\right).
\end{align}
Furthermore,
\begin{align}
  |||\mathcal Q(\infty,t)-I|||\leq 2\int_t^\infty b(\tau)\d\tau\;\exp\left(2\int_0^\infty b(\theta)\d\theta\right)\to0
\end{align}
as $t\to\infty$ implies invertibility of $\mathcal Q(\infty,t)$ for sufficiently large
$t$, while the propagation property $\mathcal Q(\infty,s)=\mathcal E_0(s-t)\mathcal Q(\infty,t) \mathcal E(t,s)$
extends this result to arbitrary $s$.

Now the statement of Theorem~\ref{thm1} follows with 
$JW_+J^*=\lim_{t\to\infty} \mathcal E_0(-t)\mathcal E(t,0)=\mathcal Q(\infty,0)$
together with the fact that $\mathcal E_0(t)$ is unitary.\hfill$\square$

\begin{rem}
  Because it was not necessary to use (A1), we can reverse the time direction in this
  statement and obtain the existence of a corresponding wave operator
  \begin{equation}
    W_-=\lim_{t\to-\infty} J^*\mathcal E_0(-t)\mathcal E(t,0)J
  \end{equation}
  when $b\in L^1(\R)$
  and construct the scattering operator $S=W_+W_-^{-1}$. This operator maps
  incoming free waves (from $t=-\infty$) to outgoing free waves (to $t=\infty$), 
  both parameterised by corresponding Cauchy data at time $t=0$,  which are 
  connected by one solution of the dissipative equation. 
\end{rem}
\begin{rem}
  Theorem~\ref{thm1} remains true if we replace \eqref{eq:CP} by
  \begin{equation}
    u''+B(t)u'+Au=0,    \tag{\ref{eq:CP}'}
  \end{equation}
  where $B=B(t)\in L^1([0,\infty),\mathcal L(H))$, cf. \cite{Wir02}. 
\end{rem}

\section{Representation of solutions}
\label{sec2}
If $b\not\in L^1[0,\infty)$ we will not apply \eqref{eq:Qdef} directly to represent
the solutions of system \eqref{eq:system} (although it is still valid). Our
approach is based on a transformation of the system \eqref{eq:system} into a 
diagonal-dominated form.

To make the calculations more transparent we use the spectral resolution
of the operator $A$ (or $\Lambda$, resp.) and reduce the analytic properties
of the operator to algebraic properties of a corresponding function. 
Following \cite[Theorem VIII.4, p.260]{RS80} there exists a unitary 
operator $\Xi:H\to L^2(X,\d\nu)$, $(X,\d\nu)$ a suitable measure space,
such that the operator $A:H\supseteq\Dom(A)\to H$ is unitarily equivalent to
\begin{equation}
  (A u)(\xi) = A(\xi) u(\xi),\qquad \xi\in X,\quad u\in\Dom(A),
\end{equation}
where $A(\xi)$ is a non-negative and $\nu$-measurable function. Furthermore,
\begin{equation}
  \Dom(A^\gamma) \simeq L^2(X, (1+A^{2\gamma}(\xi))\d\nu),\qquad \gamma\geq0.
\end{equation}
For the following the square root of $A$ is of particular importance. We denote
$\Lambda(\xi)=\sqrt{A(\xi)}$. Then equation \eqref{eq:CP} is equivalent to
\begin{equation}
  \partial_t^2u(t,\xi)+2b(t)\partial_tu(t,\xi)+\Lambda^2(\xi)u(t,\xi)=0
\end{equation}
and system \eqref{eq:system} reads as
\begin{equation}\label{eq:systemODE}
  \partial_t \begin{pmatrix} \Lambda(\xi)u\\\partial_t u \end{pmatrix}
  = \begin{pmatrix} &\Lambda(\xi)\\-\Lambda(\xi)&-2b(t) \end{pmatrix}
  \begin{pmatrix} \Lambda(\xi)u\\\partial_t u \end{pmatrix}.
\end{equation}
Estimates in the energy space $E$ correspond via $\Xi\circ J$ to estimates for this system in 
$L^2(X_+,\d\nu)\times L^2(X,\d\nu)$, where $X_+=\{\Lambda(\xi)>0\}$. 

We will solve this system following \cite{WirDiss} and \cite{Wir06}. For a 
(later to be determined) number $N>0$ we decompose the extended phase space
$[0,\infty)\times X$ in two sets, the {\em dissipative zone}
\begin{equation}
  Z_{diss}(N) := \big\{\, (t,\xi)\in[0,\infty)\times X\,\big|\, \Lambda(\xi)(1+t)< N\,\big\} 
\end{equation}
and the {\em hyperbolic zone} 
\begin{equation}
  Z_{hyp}(N) := \big\{\, (t,\xi)\in[0,\infty)\times X\,\big|\, \Lambda(\xi)(1+t)\geq  N\,\big\}.
\end{equation} 
Let further $t_\xi:X\supseteq\{\Lambda(\xi)\leq N\}\to[0,\infty)$ be implicitly defined by $\Lambda(\xi)(1+t_\xi)=N$.

In the first one we reformulate \eqref{eq:systemODE} as system of
Volterra integral equations and solve mainly by brute force, while in the
second one we apply two steps of a diagonalisation procedure to extract the 
main terms of the representation of solutions and to get sharp estimates
for the remainder terms.

\subsection{Consideration in the dissipative zone}
We will prove estimates for the fundamental solution $\mathcal E(t,s,\xi)$ to 
\eqref{eq:systemODE}, i.e. the matrix valued solution to that system with
initial condition $\mathcal E(s,s,\xi)=I\in\mathbb C^{2\times2}$. Let for this
\begin{equation}
  \lambda(t) = \exp\left(\int_0^t b(\tau)\d\tau\right)
\end{equation}
denote an auxiliary function related to the dissipation term $b(t)$. We 
distinguish different cases related to the asymptotic behaviour of $\lambda(t)$.

\medskip
\noindent{\bf Case (C1): $\overline\mu:=\limsup_{t\to\infty}tb(t)<1/2$.}

We denote by $v(t,\xi)$ and $w(t,\xi)$ the entries of one of the rows of
$\mathcal E(t,0,\xi)$. Then these functions satisfy
\begin{subequations}
\begin{align}
  v(t,\xi) &= \eta_1 + \Lambda(\xi)\int_0^t w(\tau,\xi)\d\tau,\\ 
  w(t,\xi) &= \frac1{\lambda^2(t)} \eta_2 - \Lambda(\xi)\frac1{\lambda^2(t)} \int_0^t \lambda^2(\tau)v(\tau,\xi)\d\tau,
\end{align}
\end{subequations}
where $\eta=(\eta_1,\eta_2)=(1,0)$ or $\eta=(0,1)$ for the first and second column, 
respectively.

\begin{lem}\label{lem:Zdiss1} Assume (A1) and (A2). Then
  in the case (C1) the fundamental solution $\mathcal E(t,0,\xi)$ to \eqref{eq:systemODE}
  satisfies for all $(t,\xi)\in Z_{diss}(N)$ the point-wise estimate\footnote{For a
matrix $A$ we denote by $|A|$ the matrix consisting of the absolute values of the entries of $A$.}
  \begin{equation}
    |\mathcal E(t,0,\xi)| \lesssim \frac1{\lambda^2(t)} \begin{pmatrix}
      \Lambda^{-\gamma}(\xi) & 1 \\\Lambda^{-\gamma}(\xi) & 1 \end{pmatrix}
  \end{equation}
  for any $\gamma$ such that $\Lambda^\gamma(\xi)\lambda^2(t)\lesssim1$.
\end{lem}

\begin{expl}
If we consider the case $b(t)=\mu/(1+t)$ with $\mu\in[0,\frac12)$ we obtain
$\lambda(t)\sim t^\mu$ and thus we can choose $\gamma=2\mu$ in the previous statement.
\end{expl}
\begin{expl}
In general we can calculate the constant $\gamma$ as follows: Let 
$\overline\mu=\limsup_{t\to\infty}tb(t)$ and assume $\overline\mu<\frac12$. Then the conditions of the above 
lemma are satisfied and $\lambda^2(t)\lesssim t^{2\overline\mu+\epsilon}$ for any $\epsilon>0$. Hence we choose $\gamma>2\overline\mu$. Especially, 
if $tb(t)\to0$ as $t\to\infty$ we have the above statement for any $\gamma>0$. 
Furthermore, if $b\in L^1[0,\infty)$ the choice $\gamma=0$ is possible.\footnote{Otherwise, the choice $\gamma=2\overline\mu$ 
is not always possible, as the example $b(t)=\frac1{4(e+t)}+\frac1{(e+t)\log(e+t)}$ shows.}

\end{expl}

\begin{proof}(of Lemma~\ref{lem:Zdiss1}) We start by estimating the first column, i.e. $\eta=(1,0)$.
Plugging the second integral equation into the first one yields
\begin{align}
  v(t,\xi) &=1 - \Lambda^2(\xi)\int_0^t \frac1{\lambda^2(\tau)} \int_0^\tau \lambda^2(\theta)v(\theta,\xi)\d\theta \d\tau\notag\\
  &=1 -  \Lambda^2(\xi)\int_0^t \lambda^2(\theta)v(\theta,\xi) \int_\theta^t \frac1{\lambda^2(\tau)}\d\tau \d\theta, 
\end{align}  
such that $\lambda^2(t)\Lambda^\gamma(\xi)v(t,\xi)$ satisfies an Volterra integral equation with kernel
$k(t,\theta,\xi)=-\Lambda^2(\xi)\lambda^2(t)\int_\theta^t \d\tau/{\lambda^2(\tau)}$ and source term $h(t,\xi)=\Lambda^\gamma(\xi)\lambda^2(t)\lesssim1$. We represent its solution by a Neumann series
\begin{align}
\lambda^2(t)\Lambda^\gamma(\xi)v(t,\xi)
= h(t,\xi)+\sum_{\ell=1}^\infty \int_0^t k(t,t_1,\xi)
 \cdots \int_0^{t_{k-1}} k(t_{k-1},t_k,\xi)h(t_k,\xi)\d t_k\cdots \d t_1.
\end{align}
 From (C1) we know that ${t}/{\lambda^2(t)}$ is monotone increasing for large $t$. Hence
 we conclude the kernel estimate
 \begin{align*}
  \sup_{(t,\xi)\in Z_{diss}(N)}&  \int_0^{t} \sup_{0\le\tilde t\le t_\xi} |k(\tilde t,\theta,\xi)|\d\theta
  \le \sup_{(t,\xi)\in Z_{diss}(N)}  \Lambda^2(\xi)\lambda^2(t_\xi) \int_0^t \int_\theta^{t_\xi} \frac{\d\tau}{\lambda^2(\tau)}\d\theta\\ &= \sup_{(t,\xi)\in Z_{diss}(N)}  \Lambda^2(\xi)\lambda^2(t_\xi) \int_0^{t_\xi} \frac{\tau}{\lambda^2(\tau)}\d\tau 
  \lesssim \Lambda^2(\xi)t_\xi^2\lesssim 1
 \end{align*}
which implies by standard arguments that  $|\lambda^2(t)\Lambda^\gamma(\xi)v(t,\xi)|\lesssim 1$ uniformly on $Z_{diss}(N)$. Then, the second integral equation implies
\begin{equation}
  |w(t,\xi)|\leq \frac{\Lambda^{-\gamma}(\xi)}{\lambda^2(t)} \Lambda(\xi)\int_0^t |\lambda^2(\tau)\Lambda^\gamma(\xi)v(\tau,\xi)|\d\tau \lesssim \frac{\Lambda^{-\gamma}(\xi)}{\lambda^2(t)} \Lambda(\xi)t\lesssim \frac{\Lambda^{-\gamma}(\xi)}{\lambda^2(t)}
\end{equation}
by the definition of the zone.

For the second column we use the same idea, plugging the second integral equation into the first one 
yields the new source term $\int_0^t {\Lambda(\xi)}/{\lambda^2(\tau)}\d\tau\lesssim1/{\lambda^2(t)}$, such that $\lambda^2(t)v(t,\xi)$ is uniformly bounded on $Z_{diss}(N)$.
\end{proof}

\medskip
\noindent{\bf Case (C2): $\underline\mu:=\liminf_{t\to\infty}tb(t)>1/2$.}

In this case we obtain the decay rate $1/t$ (which seems to be natural from
the point of view of effective dissipation, cf. \cite{Wir06b} and 
\cite{Yam06}).

\begin{lem}\label{lem:Zdiss2} Assume (A1) and (A2). Then
  in the case (C2)
  the fundamental solution $\mathcal E(t,0,\xi)$ to \eqref{eq:systemODE}
  satisfies for all $(t,\xi)\in Z_{diss}(N)$ the point-wise estimate
  \begin{equation}
    |\mathcal E(t,0,\xi)| \lesssim \frac1{1+t} \begin{pmatrix}
      \Lambda^{-1}(\xi) & 1 \\\Lambda^{-1}(\xi) & 1 \end{pmatrix}.
  \end{equation}
\end{lem}
\begin{proof}
We proceed in a similar way like for Lemma~\ref{lem:Zdiss1}. Plugging  the second
integral equation into the first one gives
\begin{align}
  (1+t)v(t,\xi) &= (1+t)\eta_1 +(1+t)\Lambda(\xi)\eta_2\int_0^t\frac{\d\tau}{\lambda^2(\tau)} \notag\\&\qquad
  -  \Lambda^2(\xi)\int_0^t (1+\theta)v(\theta,\xi) \int_\theta^t \frac{1+t}{1+\theta}\frac{\lambda^2(\theta)}{\lambda^2(\tau)}\d\tau \d\theta.
\end{align} 
Under condition (C2) we know that $1/\lambda^2(t)\in L^1(\R_+)$ and $\lambda^2(t)/t$
is monotone increasing for large $t$. This implies again a kernel estimate for this Volterra integral equation
\begin{align*}
  \sup_{(t,\xi)\in Z_{diss}(N)}\Lambda^2(\xi) \int_0^t &\sup_{0\le\tilde t\le t_\xi} \int_\theta^{\tilde t} \frac{1+\tilde t}{1+\theta}\frac{\lambda^2(\theta)}{\lambda^2(\tau)}\d\tau \d\theta \\ &
  \le  \sup_{(t,\xi)\in Z_{diss}(N)}\Lambda^2(\xi)  \int_0^{t_\xi} \int_0^\tau  \frac{1+t_\xi}{1+\theta}\frac{\lambda^2(\theta)}{\lambda^2(\tau)}\d\theta\d\tau
  \lesssim 1
\end{align*}
and we obtain for the first column that $|\Lambda(\xi) (1+t)v(t,\xi)|\lesssim1$ is uniformly
bounded on $Z_{diss}(N)$ while for the second one $|(1+t)v(t,\xi)|\lesssim 1$ on $Z_{diss}(N)$. 
The second integral equation yields corresponding bounds for $w(t,\xi)$ from
\begin{equation}
  \Lambda(\xi)\frac{1+t}{\lambda^2(t)}\int_0^t \frac{\lambda^2(\tau)}{1+\tau}\d\tau\lesssim (1+t) \Lambda(\xi)\lesssim 1.
\end{equation}
\end{proof}

\begin{rem}
  There remains a gap between the cases (C1) and (C2). In \cite{Wir04}
  the special case $2b(t)=\mu/(1+t)$ for $A=-\Delta$ on $L^2(\R^n)$ was studied. 
  The above exceptional value $\mu=1$ is related to the occurrence of logarithmic
  terms in the representation of solutions. 
\end{rem}

\subsection{Consideration in the hyperbolic zone}

In $Z_{hyp}(N)$ we apply two steps of diagonalisation to system 
\eqref{eq:systemODE}. In a first step we use
\begin{equation}
  M=\begin{pmatrix} i & -i\\1 &1 \end{pmatrix},
  \qquad\qquad
  M^{-1} = \frac12 \begin{pmatrix} -i & 1\\i&1 \end{pmatrix}
\end{equation}
to obtain
\begin{equation}
  M^{-1}\begin{pmatrix} & \Lambda(\xi) \\ -\Lambda(\xi)&-2b(t) \end{pmatrix} M
  = \begin{pmatrix} -i\Lambda(\xi) &\\&i\Lambda(\xi) \end{pmatrix}
  - b(t) \begin{pmatrix} 1&1\\1&1 \end{pmatrix}.
\end{equation}
In a second step we want to transform the second matrix without changing
the structure of the first one. For this we set
\begin{align}
  N^{(1)}(t,\xi) &=i 
  \begin{pmatrix} & -\frac{b(t)}{2\Lambda(\xi)} \\ \frac{b(t)}{2\Lambda(\xi)} \end{pmatrix},\\
  B^{(1)}(t,\xi) &= \partial_t N^{(1)}(t,\xi) + b(t)N^{(1)} .
\end{align}
Then we have by construction
\begin{equation}
    \begin{pmatrix} -i\Lambda(\xi) &\\&i\Lambda(\xi) \end{pmatrix} N^{(1)}(t,\xi)-N^{(1)}(t,\xi)
    \begin{pmatrix} -i\Lambda(\xi) &\\&i\Lambda(\xi) \end{pmatrix}
    = -b(t)
    \begin{pmatrix} &1\\1& \end{pmatrix},
\end{equation}
such that with  $N_1(t,\xi) = I-N^{(1)}(t,\xi)$ the operator identity
\begin{multline*}
  \left(\partial_t- \begin{pmatrix} -i\Lambda(\xi) &\\&i\Lambda(\xi) \end{pmatrix}+
    b(t) \begin{pmatrix} 1&1\\1&1 \end{pmatrix}\right) N_1(t,\xi) \\
  = \partial_t - N_1(t,\xi)\begin{pmatrix} -i\Lambda(\xi) &\\&i\Lambda(\xi) \end{pmatrix}+
  b(t)N_1(t,\xi)-B^{(1)}(t,\xi)
\end{multline*}
holds. If the zone constant $N$ is chosen sufficiently large, we have that 
$\det N_1 (t,\xi)=1-{b^2(t)}/{(4\Lambda^2(\xi))}\geq c>0$ is uniformly bounded away from zero
on $Z_{hyp}(N)$. Then $N_1^{-1}(t,\xi)$ exists and $N_1(t,\xi)$ and $N_1^{-1}(t,\xi)$ are
both  uniformly bounded on $Z_{hyp}(N)$. 

Setting $R_1(t,\xi)=-N_1^{-1}(t,\xi) B^{(1)}(t,\xi)$ we obtain the operator identity
\begin{multline}
  N_1^{-1}(t,\xi)\left(\partial_t- \begin{pmatrix} -i\Lambda(\xi) &\\&i\Lambda(\xi) \end{pmatrix}+
    b(t) \begin{pmatrix} 1&1\\1&1 \end{pmatrix}\right) N_1(t,\xi) \\
  =\partial_t - \begin{pmatrix} -i\Lambda(\xi) &\\&i\Lambda(\xi) \end{pmatrix}+
  b(t)I + R_1(t,\xi)
\end{multline}
with remainder term $R_1(t,\xi)$ subject to the point-wise estimate
\begin{equation}
  \|R_1(t,\xi)\|\lesssim\frac1{\Lambda(\xi)(1+t)^2}
\end{equation}
following directly from Assumption (A2). Note, that we did not use Assumption (A1) for the treatment
of the hyperbolic zone.

We are now in a position to derive the main result of this section.

\begin{lem}\label{lem:ZhypRep} Assume (A2).
  Then the fundamental solution $\mathcal E(t,s,\xi)$ of \eqref{eq:systemODE}
  can be represented as
  \begin{equation}
    \mathcal E(t,s,\xi) =\frac{\lambda(s)}{\lambda(t)} M^{-1}N_1^{-1}(t,\xi)
       \widetilde{\mathcal E_0}(t,s,\xi) \mathcal Q_1(t,s,\xi) N_1(s,\xi)M
  \end{equation}
  for $t\geq s$ and $(s,\xi)\in Z_{hyp}(N)$, where
  \begin{itemize}
  \item the function $\lambda(t)=\exp\left(\int_0^t b(\tau)\d\tau\right)$ describes the
    main influence of the dissipation $b=b(t)$,
  \item the matrices $N_1(t,\xi)$ and $N_1^{-1}(t,\xi)$ are uniformly bounded
    on $Z_{hyp}(N)$ tending on $\{\Lambda(\xi)\geq\epsilon\}\subseteq X$ uniformly to the identity matrix $I$,
  \item the matrix $\widetilde{\mathcal E_0}(t,s,\xi)$ is given by
    \begin{equation}
      \widetilde{\mathcal E_0}(t,s,\xi)=
      \begin{pmatrix} e^{-i\Lambda(\xi)(t-s)}&\\&   e^{i\Lambda(\xi)(t-s)} \end{pmatrix},
    \end{equation}
  \item and the matrix $\mathcal Q_1(t,s,\xi)$ is uniformly bounded and invertible 
    on $Z_{hyp}(N)$ tending on $\{\Lambda(\xi)\geq\epsilon\}\subseteq X$ uniformly to the invertible
    matrix $\mathcal Q_1(\infty,s,\xi)$.
  \end{itemize}
\end{lem}

\begin{rem}
  The free propagator can be represented as 
  $\mathcal E_0(t-s,\xi)=M^{-1}\widetilde{\mathcal E_0}(t,s,\xi)M$. Thus, 
  Lemma~\ref{lem:ZhypRep} gives a relation between the free propagator
  $\mathcal E_0(t-s,\xi)$ and $\mathcal E(t,s,\xi)$.
\end{rem}
\begin{rem}
  If the operator $A$ is boundedly invertible, the function $\Lambda(\xi)$ is 
  bounded from below and hence the matrices $N_1(t,\xi)$ and $\mathcal Q_1(t,\xi)$
  converge uniformly on $X$ as $t\to\infty$. Furthermore, the dissipative zone can be skipped in
  this case (because it is contained in a finite time strip).  Thus Lemma~\ref{lem:ZhypRep} describes all
  essential properties of the fundamental solution in this case.
\end{rem}

\begin{proof}(of Lemma~\ref{lem:ZhypRep})
The construction of the representation of solutions will be done in two
steps. First, note that
\begin{equation}
 \Tilde{\Tilde{\mathcal E}}_0(t,s,\xi) = \frac{\lambda(s)}{\lambda(t)}  \widetilde {\mathcal E_0}(t,s,\xi)
\end{equation}
is the fundamental solution to the diagonal main part,
\begin{equation}
 \partial_t \Tilde{\Tilde{\mathcal E}}_0(t,s,\xi)
 = \begin{pmatrix} -i\Lambda(\xi)&\\&i\Lambda(\xi)\end{pmatrix} \Tilde{\Tilde{\mathcal E}}_0(t,s,\xi)
 + b(t) \Tilde{\Tilde{\mathcal E}}_0(t,s,\xi),\qquad  \Tilde{\Tilde{\mathcal E}}_0(s,s,\xi)=I. 
 \end{equation} 
Thus making the ansatz $\Tilde{\Tilde{\mathcal E}}_0(t,s,\xi)\mathcal Q_1(t,s,\xi)$ for the fundamental solution to the 
transformed operator $\partial_t -\mathrm{diag}\left( -i\Lambda(\xi),i\Lambda(\xi)\right)+ b(t)I+R_1(t,\xi)$, we obtain a system 
for $\mathcal Q_1(t,s,\xi)$,
\begin{equation}
  \partial_t \mathcal Q_1(t,s,\xi)=\mathcal R(t,s,\xi)\mathcal Q_1(t,s,\xi),\qquad \mathcal Q_1(s,s,\xi)=I
\end{equation}
with coefficient matrix $\mathcal R(t,s,\xi)= \widetilde{\mathcal E_0}(s,t,\xi)R_1(t,\xi)\widetilde{\mathcal E_0}(t,s,\xi)$.
Using that $\widetilde{\mathcal E_0}(t,s,\xi)$ is unitary, we see that $\mathcal R(t,s,\xi)$ satisfies the
same estimates as $R_1(t,\xi)$. 

This allows us to estimate in a second step
the solution $\mathcal Q_1(t,s,\xi)$ directly from the representation as Peano-Baker series
\begin{equation}\label{eq:Q1-ser}
  \mathcal Q_1(t,s,\xi)=I+\sum_{k=1}^\infty \int_s^t \mathcal R(t_1,s,\xi)\int_s^{t_1}\mathcal R(t_2,s,\xi)\cdots\int_s^{t_{k-1}}\mathcal R(t_k,s,\xi)\d t_k\cdots \d t_1.
\end{equation}
Thus,
\begin{equation}
  \|\mathcal Q_1(t,s,\xi)\| \leq \exp\left(\int_s^t \|R_1(\tau,\xi)\|\d\tau\right) \leq \exp\left(\frac c{\Lambda(\xi)} \int_s^t \frac{\d\tau}{(1+\tau)^2}\right) \lesssim 1
\end{equation}
uniformly in $t\geq s$ and $(s,\xi)\in Z_{hyp}(N)$. Note that the inverse matrix satisfies 
a similar differential equation and therefore an analogue to series \eqref{eq:Q1-ser}. 

It remains to show that $\mathcal Q_1(t,s,\xi)$ converges uniformly on
$\{\,\Lambda(\xi)\geq\epsilon\,\}\subseteq X$ as $t\to\infty$. This follows by the Cauchy criterion applied to the series \eqref{eq:Q1-ser} or by the estimate
\begin{multline}\label{eq:2.30}
  \|\mathcal Q_1(\infty,s,\xi)-\mathcal Q_1(t,s,\xi)\|\leq \int_t^\infty \|R_1(\tau,\xi)\| \d\tau  \exp\left(\int_s^t \|R_1(\tau,\xi)\|\d\tau\right)\\
  \lesssim \frac1{\Lambda(\xi)} \int_t^\infty \frac{\d\tau}{(1+\tau)^2} \lesssim \frac1{\epsilon(1+t)}. 
\end{multline}
The theorem is proven.
\end{proof}

\subsection{Combination of results}
We combine the previously constructed representations of $\mathcal E(t,s,\xi)$
to get a representation for $\mathcal E(t,0,\xi)$ and the corresponding operator
$\mathcal E(t,0)$. We will obtain $J^{*}\mathcal E(t,0)J\in\mathcal L(E)$
and
\begin{equation}\label{eq:Eest}
  |||J^{*}\mathcal E(t,0)J|||_{\tilde E \to E}\lesssim \frac1{\lambda(t)}
\end{equation}
by the right choice of a smaller space $\tilde E\subset E$. Again $J$ denotes the canonic
identification $E\simeq (H\ominus\Ker\Lambda)\times H$.

If $\Lambda(\xi)>N$ for all $\xi\in X$, we have to consider only the hyperbolic zone and 
$\mathcal E(t,0,\xi)$ is given by Lemma~\ref{lem:ZhypRep}, hence \eqref{eq:Eest} holds true
with $\tilde E=E$. We focus on the remaining case and consider $\xi$ with $\Lambda(\xi)\le N$ here. There it holds
\begin{equation}
  \mathcal E(t,0,\xi)=\mathcal E(t,t_\xi,\xi)\mathcal E(t_\xi,0,\xi)
\end{equation}
and we obtain corresponding estimates by Lemmata~\ref{lem:Zdiss1} to \ref{lem:ZhypRep}. We distinguish some cases.

\medskip\noindent
{\bf Case (C1):} Using Lemma~\ref{lem:Zdiss1} and \ref{lem:ZhypRep} we conclude
\begin{equation}\label{eq:C1_Eest}
  |\mathcal E(t,0,\xi)|\lesssim\frac{\lambda(t_\xi)}{\lambda(t)}
   \frac1{\lambda^2(t_\xi)} \begin{pmatrix}
      \Lambda^{-\gamma}(\xi) & 1 \\\Lambda^{-\gamma}(\xi) & 1 \end{pmatrix}
    \lesssim \frac1{\lambda(t)} \begin{pmatrix}
      \Lambda^{-\gamma}(\xi) & 1 \\\Lambda^{-\gamma}(\xi) & 1 \end{pmatrix}
\end{equation}
for $\Lambda(\xi)\le N$, $t\ge t_\xi$ and  
$\gamma$ with $\Lambda^\gamma(\xi)\lambda^2(t_\xi)\lesssim 1$.  For $t\le t_\xi$ the same estimate
follows from Lemma~\ref{lem:Zdiss1}.

To simplify notation
we introduce $[\Lambda(\xi)]=\min(\Lambda(\xi),N)$. If we choose $\gamma=1$, we need for the data 
an estimate of the form $\|[\Lambda(\xi)]^{-1}\Lambda(\xi)u_1(\xi)\|_2\sim\|\langle\Lambda(\xi)\rangle u_1(\xi)\|_2=\|u_1\|_{\Dom(\Lambda)}$.  Thus we obtain  that
the operator $J^*\mathcal E(t,0)J$ restricted to $\Dom(\Lambda)\times H\to E$ satisfies the
bound \eqref{eq:Eest}.

\medskip\noindent
{\bf Case (C2):} Lemma~\ref{lem:Zdiss2} implies only the weaker decay rate 
$1/(1+t)$. This can be compensated by assumptions on the data, which are
valid near $\{\Lambda(\xi)=0\}\subseteq X$. Using
\begin{equation}
  \frac{\Lambda^\gamma(\xi)\lambda(t_\xi)}{1+t_\xi}\sim\frac{\lambda(t_\xi)}{(1+t_\xi)^{\gamma+1}}\lesssim1 
\end{equation}
for $\gamma>\overline\mu-1$, $\overline\mu=\limsup_{t\to\infty}tb(t)$, we conclude that
\begin{equation}\label{eq:C2_Eest}
 |\mathcal E(t,0,\xi)|
 \lesssim \frac{\lambda(t_\xi)}{\lambda(t)} \frac1{1+t_\xi} \begin{pmatrix}
       \Lambda^{-1}(\xi) & 1 \\\Lambda^{-1}(\xi) & 1
 \end{pmatrix}
 \lesssim\frac1{\lambda(t)} \begin{pmatrix}
      \Lambda^{-1}(\xi) & 1 \\ \Lambda^{-1}(\xi) & 1 \end{pmatrix} \Lambda^{-\gamma}(\xi)
\end{equation}
for all $\Lambda(\xi)\le N$ and $t\ge t_\xi$. For $t\le t_\xi$ Lemma~\ref{lem:Zdiss2} implies the same bound using $\Lambda^\gamma(\xi)\lambda(t)/(1+t)\lesssim \lambda(t)/(1+t)^{\gamma+1}\lesssim 1$.

Thus if we define for $\gamma\geq0$ the 
{\em modified energy space} $E^{(\gamma)}$ to be the closure of 
$(\Dom(\Lambda)\cap R(\Lambda^{\gamma})\ominus\Ker\Lambda)\times (R(\Lambda^{\gamma})\ominus\Ker\Lambda)$ with respect to the norm
\begin{equation}\label{eq:EgammaDef}
  \| (u_1,u_2) \|_{E^{(\gamma)}}^2 = \left\| [\Lambda(\xi)]^{-\gamma-1}\Lambda(\xi) u_1(\xi) \right\|^2_{L^2(X,\d\nu)} 
  +  \| [\Lambda(\xi)]^{-\gamma} u_2(\xi) \|_{L^2(X,\d\nu)}^2,
\end{equation}
the operator $J^*\mathcal E(t,0)J$ restricted to $E^{(\gamma)}\to E$ satisfies the
bound \eqref{eq:Eest}. 

\begin{rem}
  Using $[\Lambda(\xi)]^{-1}\Lambda(\xi)\sim\langle\Lambda(\xi)\rangle$ we obtain especially $E^{(0)}=(\Dom(\Lambda)\times H)\ominus(\Ker\Lambda)^2$. So in the case (C1) we may use $E^{(0)}$ for the
  data, while in (C2) we use $E^{(\gamma)}$ with $\gamma=\max(\overline\mu^+-1,0)$, where
  $\overline\mu^+$ denotes any number larger than $\overline\mu$.
\end{rem}

\begin{thm}\label{thm:EnEst}
  Let $u$ be the solution to \eqref{eq:CP} for data $(u_1,u_2)\in E^{(\gamma)}$
   for $\gamma=\max(\overline\mu^+-1,0)$. Then the estimate
  \begin{equation}
    \|(u,u')\|_E \lesssim \frac1{\lambda(t)}  \|(u_1,u_2)\|_{E^{(\gamma)}}
  \end{equation}
  is valid under assumptions (A1), (A2) in the cases (C1) and (C2).
\end{thm}

\begin{rem}
  We used conditions on the data related to $\{ \Lambda(\xi)=0\}$. In principle there
  arise two different cases. On the one hand, if $\nu\{\Lambda(\xi)=0\}>0$ 
  it follows that $0\in\sigma_P(A)$ and $0$ is an eigenvalue of $A$. In this case
  we cut out the null-space of $A$ and $R(A)$ is dense in the 
  ortho-complement  of $\Ker A$. Hence $E^{(\gamma)}$ is dense in $E^{(0)}=
  (\Dom(\Lambda)\ominus\Ker\Lambda)\times(H\ominus\Ker\Lambda)$, which is dense in the ortho-complement of $\{0\}\times\Ker\Lambda$ in $E$.

  If $\nu\{\Lambda(\xi)=0\}=0$ is a null-set, $R(A)$ is dense in $H$ and the result
  of Theorem~\ref{thm:EnEst} holds on a dense subset of $E^{(0)}=\Dom(\Lambda)\times H$.
\end{rem}

\section{Scattering type theorems}
\label{sec3}

The main purpose of this note is to explain in which sense
Theorem~\ref{thm:EnEst} is sharp. This sharpness is formulated as
a generalisation of Theorem~\ref{thm1} with an additional energy decay
function.

\begin{thm}\label{thm:3.1}
  Under the assumptions of Theorem~\ref{thm:EnEst} there exists 
  a bounded operator $W_+:E^{(\gamma)}\to E$, $\gamma=\max(\overline\mu^+-1,0)$
  such that for Cauchy data
  $(u_1,u_2)\in E^{(\gamma)}$ of \eqref{eq:CP} and associated data $(v_1,v_2)=W_+(u_1,u_2)$
  to \eqref{eq:CPfree} the corresponding solutions $u=u(t)$ and $v=v(t)$
  satisfy 
  \begin{equation}
    \| \lambda(t)(u,u')-(v,v')\|_E \to 0
  \end{equation}
  as $t\to\infty$. Furthermore, the operator $W_+$ satisfies $\Ker W_+ = \{0\}$.
\end{thm}

\begin{proof}
  The proof is based on an explicit representation of the (modified) inverse wave
  operator $W_+$. From Lemma~\ref{lem:ZhypRep} we know that the limit
  $\mathcal Q_1(\infty,t_\xi,\xi)=\lim_{t\to\infty}\mathcal Q_1(t,t_\xi,\xi)$ exists uniformly in 
  $\{ \Lambda(\xi)\geq\epsilon\}$ for any $\epsilon>0$. Hence, if we consider 
  $\mathcal E_0(-t,\xi)\mathcal E(t,0,\xi)$ on $\{ \Lambda(\xi)\geq\epsilon\}$ we obtain
  \begin{align*}
     &\lim_{t\to\infty} \lambda(t)\mathcal E_0(-t,\xi)\mathcal E(t,0,\xi) \\
     =& \lim_{t\to\infty} \lambda(t_\xi) \mathcal E_0(-t,\xi) M^{-1}N_1^{-1}(t,\xi)
     \widetilde{\mathcal E_0}(t,t_\xi,\xi)\mathcal Q_1(t,t_\xi,\xi)N_1(t_\xi,\xi)M
     \mathcal E(t_\xi,0,\xi)\\
     =& \lim_{t\to\infty} \lambda(t_\xi) \mathcal E_0(-t,\xi) \mathcal E_0(t-t_\xi,\xi)
     M^{-1}\mathcal Q_1(t,t_\xi,\xi)N_1(t_\xi,\xi)M
     \mathcal E(t_\xi,0,\xi)\\
     =& \lambda(t_\xi) \mathcal E_0(-t_\xi,\xi) M^{-1}\mathcal Q_1(\infty,t_\xi,\xi)N_1(t_\xi,\xi)M  
     \mathcal E(t_\xi,0,\xi)
  \end{align*}
  using $N_1^{-1}(t,\xi)\to I$ uniformly on $\{ \Lambda(\xi)\geq\epsilon\}$.
  Denoting this limit by $\tilde W_+(\xi)$ we see
  \begin{equation}
    |\tilde W_+(\xi)|\lesssim \lambda(t_\xi)|\mathcal E(t_\xi,0,\xi)|.
  \end{equation}
  Now the estimates \eqref{eq:C1_Eest} and \eqref{eq:C2_Eest} on 
  $\mathcal E(t_\xi,0,\xi)$ imply that $\tilde W_+(\xi)$ defines on 
  $\{\Lambda(\xi)>0\}=X\setminus\{\Lambda(\xi)=0\}$ a multiplication operator with
  estimate
  \begin{equation}
    |\tilde W_+(\xi)|\lesssim \begin{pmatrix}
      \Lambda^{-1}(\xi) & 1 \\\Lambda^{-1}(\xi) & 1 \end{pmatrix} [\Lambda(\xi)]^{-\gamma}
  \end{equation}
  for $\gamma=\max(\overline\mu^+-1,0)$. This implies that the corresponding operator
  $W_+=J^*\tilde W_+J$ resticted to $E^{(\gamma)}\to E$ is bounded. Furthermore, we obtain from the fact that the free
  propagator $\mathcal E_0(t)$ is unitary for all $t$
  \begin{multline}
    \|(u_1,u_2)\|_{E^{(\gamma)}} \gtrsim \| \lambda(t)\mathcal E(t,0)(u_1,u_2) - \mathcal E_0(t) W_+ (u_1,u_2)\|_E\\
    = \|\lambda(t)\mathcal E_0(-t)\mathcal E(t,0)(u_1,u_2) - W_+(u_1,u_2)\|_E \to0
  \end{multline}
  for all $(u_1,u_2)$ with $\mathrm{supp}_\xi (u_1,u_2)\subseteq\{\Lambda(\xi)>0\}$. But this
  is a dense subset of $E^{(\gamma)}$ and by Banach-Steinhaus
  theorem this gives convergence for all $(u_1,u_2)$. 

  On $\{\Lambda(\xi)>0\}$ the inverse wave operator $W_+$ is represented as product of invertible   
  matrices. Since $\{\Lambda(\xi)=0\}$ is excluded by the definition of $E^{(\gamma)}$, the null-space 
  $\Ker W_+$ of $W_+:E^{(\gamma)}\to E$ is trivial  and the theorem  is proven.
\end{proof}

If we combine the previous theorem with the conservation of abstract energy for
the free problem \eqref{eq:CPfree}, $\|(v,v')\|_E=\|(v_1,v_2)\|_E$, we immediately obtain

\begin{cor}
  For all data $(u_1,u_2)\in E^{(\gamma)}$ with $\gamma=\max\{\overline\mu^+-1,0\}$
  the abstract energy of the solution to \eqref{eq:CP} satisfies the two-sided estimate
  \begin{equation}\label{eq:3.5}
    \|(u,u')\|_E \sim \frac1{\lambda(t)}
  \end{equation}
  for all $t\geq0$.
\end{cor}

Note, that the definition of $E^{(\gamma)}$ (formula~\eqref{eq:EgammaDef})
implies $E^{(\gamma)}\perp\Ker A$, such that corresponding eigenfunctions of $A$ are excluded 
in the previous corollary. If the data belong to $\Ker A$, the energy decays like $1/\lambda^2(t)$.

If the operator $A$ has a bounded inverse, i.e. $\Lambda(\xi)\geq c_0>0$ uniformly on $X$, the dissipative
zone $Z_{diss}(N)$ can be shrunken to $Z_{diss}(N)\cap\{\Lambda(\xi)\geq c_0\}$, such that $1+t\leq Nc_0^{-1}$ holds uniformly.
In this case, the estimate of $\mathcal E(t,s,\xi)$ in the dissipative zone can be replaced by the trivial
estimate $\|\mathcal E(t,s,\xi)\|\lesssim1$ following directly from system \eqref{eq:systemODE} using the boundedness of
$|b(t)|$ and $\Lambda(\xi)$. Thus, in this case we can revert the time variable and consider the limits as $t\to-\infty$ under (A2)  as
well. Furthermore, $E=E^{(\gamma)}=\Dom(\Lambda)\times H$ for all choices of $\gamma$ and Theorem~\ref{thm:3.1} follows directly
from Lemma~\ref{lem:ZhypRep} without relying on Banach-Steinhaus theorem. 

\begin{thm}\label{thm:3.3}
  Assume $(Aw,w)\geq c_0\|w\|$ for all $w\in H$. Assume further that the coefficient $b=b(t)$ satisfies (A2). 
  Then there exist bounded and invertible operators $W_+:E\to E$ and $W_-:E\to E$, such that for
  all Cauchy data $(u_1,u_2)\in E$ and corresponding data $(v_1^\pm,v_2^\pm)=W_\pm(u_1,u_2)$
  the estimate
  \begin{equation}
    \| \lambda(t)(u,u')-(v,v')\|_E \le C \langle t\rangle^{-1} \|(u_1,u_2)\|_E
  \end{equation}
  holds true. The constant $C$ depends on the coefficient $b(t)$ and $c_0$.
\end{thm}
\begin{proof}
  The statement follows directly from Lemma 2.3 (for positive times and if we replace $b(t)$ by $-b(t)$ 
  also for negative times) if we recall the following additional estimates (written for $t\ge0$) 
  \begin{equation}
     \|N_1^{-1}(t,\xi)-I\| + \|\mathcal Q_1(t,s,\xi)-\mathcal Q_1(\infty,s,\xi)\|
   \lesssim \frac1{\Lambda(\xi)(1+t)} \lesssim \frac1{1+t},
   \end{equation}
  The first one follows from the properties of $N^{(1)}(t,\xi)$, the last one restates \eqref{eq:2.30}. 
  
  The operator $W_+$ is defined via the limit $\tilde W_+(\xi)$ as in the proof of Theorem~\ref{thm:3.1},
  but the lower bound on $A$ allows to skip the dissipative zone and to use $t_\xi=0$ instead.
  Then the above estimates imply for the difference
  \begin{align*}
     &\|\lambda(t) \mathcal E(t,0,\xi) - \mathcal E_0(t) W_+(\xi)\| \\
     & = \| \underbrace{M^{-1} N_1^{-1}(t,\xi)M}_{=I+\mathcal O(t^{-1})}
     \mathcal E_0(t,\xi)\underbrace{M^{-1}\mathcal Q_1(t,0,\xi)N_1(0,\xi)M}_{=(\ast)+\mathcal O(t^{-1})}\\ 
     &\qquad \qquad \qquad \qquad  - \mathcal E_0(t,\xi) \underbrace{M^{-1} \mathcal Q_1(\infty,0,\xi)N_1(0,\xi)M}_{=(\ast)}\|\\
      &\lesssim \frac1{1+t}
  \end{align*}
  and the desired statement follows.
\end{proof}

\begin{rem}
In contrast to Theorem~\ref{thm:3.1} we can use Theorem~\ref{thm:3.3} to obtain two-sided norm estimates for $J^*\mathcal E(t,0)J$. It follows that 
  \begin{equation}
    |||J^*\mathcal E(t,0)J|||_{E\to E} \sim \frac1{\lambda(t)}
  \end{equation}
  for all $t$.
\end{rem}

\section{Applications}
\label{sec4}

We will sketch several applications to the previously given abstract results.
The list is not complete, but it is intended to give an overview over the 
different cases and types of results.

\begin{expl}
  As first example we consider $A=-\Delta$ on $L^2(\R^n)$. Then the abstract 
  equation \eqref{eq:CP} reads as wave equation with time-dependent 
  dissipation. This case was treated in \cite{WirDiss} and \cite{Wir06}
  by similar techniques. The abstract spaces reduce to Sobolev scale,
  we have $\Dom(\Lambda)=H^1(\R^n)$ and $E=\dot H^1(\R^n)\times L^2(\R^n)$. 

  In this special case $\Ker \Delta=\{0\}$. Thus, as special result, we obtain that
  for arbitrary data $u_1\in H^1(\R^n)$ and $u_2\in L^2(\R^n)$ the hyperbolic energy
  of the corresponding solution $u=u(t,x)$ satisfies
  \begin{equation}\label{eq:TwoSidedEnEst}
    \|\nabla u(t,\cdot)\|_2^2+\|u_t(t,\cdot)\|_2^2 \sim \frac1{\lambda^2(t)}
  \end{equation}
  as two-sided energy estimate.

  The same result is true for $A=-\Delta$ the Dirichlet- or Neumann-Laplacian on 
  $L^2(\Omega)$ for an exterior domain $\Omega\subseteq\R^n$ with smooth boundary.
\end{expl}

\begin{expl}
  As second example we consider the Klein-Gordon type equation, where
  $A=-\Delta+1$ on $L^2(\R^n)$. Then $\sigma(A)=[1,\infty)$ and $A$ is boundedly invertible.
  Furthermore, $\Dom(\Lambda)=H^1(\R)$ and $E=H^1(\R^n)\times L^2(\R^n)$ and $E^{(\gamma)}=E$ for all
  $\gamma\geq0$.
  
  In this case it is a simple consequence of the representation of 
  $\mathcal Q_1(t,s,\xi)$ that $W_+\in\mathcal L(E)$ is invertible in 
  $\mathcal L(E)$. Especially \eqref{eq:TwoSidedEnEst} is valid.

  The same result holds true in the case of $A=-\Delta$ the Dirichlet-Laplacian
  on $L^2(\Omega)$, $\Omega$ an interior domain with smooth boundary.
\end{expl}

\begin{expl}
  Contrary to this for $A=-\Delta$ the Neumann-Laplacian on $L^2(\Omega)$, $\Omega$
  bounded domain with smooth boundary, $A$ has the non-trivial null-space
  $\Ker A=\mathbb C$ consisting of constant functions. In this case
  the two sided estimate \eqref{eq:TwoSidedEnEst} does not hold. Solving
  the corresponding ordinary differential equation we see that
  solutions to constant data behave like $1/ \lambda^2(t)$ or $b(t)/\lambda^2(t)$.
\end{expl}

\begin{expl}
  Let $(a_{ij}(x))_{i,j=1,\ldots n}\in L^\infty(\R^n,\mathbb C^{n\times n})$ be a positive and self-adjoint 
  matrix. Then we can consider the second order elliptic operator
  \begin{equation}
    A = -\sum_{i,j=1,\ldots n} \partial_i a_{ij}(x) \partial_j 
  \end{equation}
  on $L^2(\R^n)$ with domain $\Dom(A)=H^2(\R^n)$. The corresponding abstract
  energy space has norm
  \begin{equation}
    \| (u_1,u_2) \|_E^2  = \int_{\R^n} \left(|T(x)\nabla u|^2 + |u_t|^2\right) \d x  
  \end{equation}
  where $T(x)=\sqrt{(a_{ij}(x))}$ in the sense of positive matrices.
\end{expl}

\begin{expl}
  The approach is not restricted to second order operators. We can also consider
  for example the damped plate equation with $A=\Delta^2$. In this case the abstract energy
  space is $E=\dot H^2(\R^n)\times L^2(\R^n)$.
\end{expl}

\bigskip\noindent
{\bf Acknowledgements.} The author thanks the referee for his many valuable comments to improve
this paper.

\end{document}